\documentclass[]{article} 
\usepackage{spie}
\usepackage{amsmath}
\usepackage{amssymb}

\title{Classifying Characteristic  functions giving Weyl-Heisenberg
Frames}

\author{ P.G. Casazza \supit{a} and M.C. Lammers\supit{b} 
\skiplinehalf 
\supit{a} Dept. of Mathematics, University of Missouri, Columbia, MO
65203 
\\
\supit{b}  Dept. of Mathematics, University of South Carolina, 
Columbia, SC 29208 }

\authorinfo{ E-mail:pete@math.missouri.edu, lammers@math.sc.edu }

\newcommand{\lb}{\left \langle }
\newcommand{\rb}{\right \rangle}

\newcommand{\no}{\nonumber}
\newcommand{\be}{\begin{equation}}
\newcommand{\ee}{\end{equation}}

\newcommand{\bea}{\begin{eqnarray}}
\newcommand{\eea}{\end{eqnarray}}
\newcommand{\bean}{\begin{eqnarray*}}
\newcommand{\eean}{\end{eqnarray*}}
\newcommand{\ba}{\begin{array}}
\newcommand{\ea}{\end{array}}
\newcommand{\ol}{\overline}

\newcommand{\fkb}{{\frac{k}{b}}}
\newcommand{\fjb}{{\frac{j}{b}}}

\newcommand{\sfb} {{\sqrt{\frac{1}{b}} }}
\newcommand{\fb}{{ \frac{1}{b}}} 
\newcommand{\lr}{{ L^2 (\mathbb  R)}}

\newcommand{\mr}{{\mathbb  R}}
\newcommand{\mz}{{\mathbb  Z}}

\newcommand{\bd}{\begin{displaystyle}}
\newcommand{\ed}{\end{displaystyle}}

\newtheorem{Thm}{Theorem}[section]

\newtheorem{Prop}[Thm]{Proposition}
\newtheorem{Def}[Thm]{Definition}
\newtheorem{Rmk}[Thm]{Remark}

\newtheorem{Ex}[Thm]{Example}
\newtheorem{Prob}[Thm]{Problem}

\def \Bbb{\mathbb}

\begin{document} 
 \maketitle

\begin{abstract}
We examine the question of which characteristic functions
yield Weyl-Heisenberg frames for various values of the
parameters.  We also give numerous applications of frames of
characteristic functions to the general case $(g,a,b)$. 
\end{abstract}


\section{Introduction}
\label{sect:intro} 

In 1952 Duffin And Schaeffer \cite{DS} introduced the notion of
a frame for a Hilbert space.

\begin{Def}
A sequence $(f_{n})$ in a Hilbert space $H$ is a {\bf frame} for
$H$ if there are constants $0<A,B$ satisfying
$$
A\|f\|^{2}\le \sum_{n}|\langle f,f_{n}\rangle |^{2} \le
B\|f\|^{2},\ \ \text{for all}\ \ f\in H.
$$
\end{Def}

The numbers $A,B$ are called  {\bf lower} (resp. {\bf upper}) 
{\bf frame bounds for the frame}.  If $A=B$ we call this a {\bf tight
frame} and if $A=B=1$ we call it a {\bf normalized tight frame}.  If
$(f_{n})$ does not span $H$ but is a frame for its closed linear span,
we call it a {\bf frame sequence}.  It is clear from the definition
that an orthogonal projection takes a frame to a frame sequence with
the same frame bounds. 

If $(f_{n})$ is a sequence of elements of an infinite dimensional
Hilbert space $H$ and $(e_{n})$ is an orthonormal basis for $H$,
we define the {\bf preframe  operator} $T:H\rightarrow H$ by:
$Te_{n} = f_{n}$.  It follows that for any $f\in H$, $T^* f =
\sum_{n}\langle f,f_{n}\rangle e_{n}$.  Hence, $(f_{n})$ is a 
frame if and only if $T^* $ is an isomorphism (called the
{\bf frame transform}) and in this case $S=TT^*$ is an invertible
operator on $H$ called the {\bf frame operator}.  The frame operator
is a positive, self-adjoint invertible operator on $H$ satisfying:
$$
Sf = \sum_{n}\langle f,f_{n}\rangle f_{n}.
$$

A bounded unconditional basis for $H$ is called a {\bf Riesz basis}
(or a {\bf Riesz basic sequence} if it is a Riesz basis for its closed
linear span in $H$).

The frames commonly used in signal processing are the Weyl-Heisenberg 
(or Gabor) frames.  If $g\in \lr$ and $0<a,b\in \mathbb R$ we define
$$
\text{Translation by a }\ \ \ T_{a}(g)(t) = g(t-a)\ \ \ \ \ 
\text{ Modulation by b }\ \ \ E_{b}g(t) = e^{2{\pi}imbt}g(t).
$$
We say that $(g,a,b)$ generates (or {\bf is}) a {\bf Weyl-Heisenberg
frame} 
({\bf WH-frame} for short) for $\lr$ if
$(E_{mb}T_{na}g)_{m,n\in \mathbb Z}$ is a frame for $\lr$.  If this
family has a finite upper frame bound we call $g$ a {\bf preframe 
function}.  The family of preframe  functions is denoted {\bf PF}.
 
There are several known restrictions on the $g,a,b$ in order that
$(g,a,b)$ form a WH-frame which we summarize below.  These results
are due to various authors and may be found in Heil and Walnut
\cite{HW} or Casazza \cite{C}.
To simplify the notation, for all $k\in \mathbb Z$ we let
$$
G_{k}(t) = \sum_{n\in \mathbb Z}g(t-na)\overline{g(t-na-\fkb )}.
$$

\begin{Prop} \label{P1}
Let $g\in${\bf PF} and $0<a,b$.  

(1)  If $(g,a,b)$ generates a WH-frame then $ab\le 1$.

(2)  $(g,a,b)$ generates a WH-frame if and only if $(g,\fb,\frac{1}{a})$
generates a Riesz basic sequence.  Hence, if $ab=1$, then
$(g,a,b)$ is a WH-frame if and only if it is a Riesz basis of $\lr$.

(3)  If $(g,a,b)$ is a WH-frame with frame bounds $A,B$ then
$bA\le G_{0}(t)\le bB$ a.e.

(4)  If $ab\le 1$ and supp g $\subset [0,\fb]$, then $(g,a,b)$ is a
WH-frame with frame bounds $A,B$ if and only if $bA\le G_{0}(t)\le bB$
a.e.

(5)  If $a=b=1$, and for a.e. $0\le t\le 1$ we have that $(g(t-n))_{n\in
\mathbb Z}$ has at most one non-zero term, then $(g,1,1)$ is a WH-frame
with
frame bounds $A,B$ if and only if $A\le G_{0}(t)\le B$ a.e.

(6) If $(g,a,b)$ is in PF with upper frame bound B then
$\sum|g(x-\fkb)|^2 \le
B$.
\end{Prop}

There is a necessary condition for having a WH-frame due to Casazza and
Christensen \cite{CC1}. 

\begin{Thm}[{\bf CC-Condition}]\label{pfcc}
Let $g\in \lr$, $a,b>0$ and assume that\\
(1)We have  
\[
A=:\inf_{t\in [0,a]}
\left [G_0(t) - \sum_{k\not= 0}|G_{k}(t)|\right ] > 0, \]
(2)  We have,
$$
B=:\sup_{t\in [0,a]}\sum_{k\in  \mz }|G_{k}(t)|< \infty.
$$
Then $(E_{mb}T_{na}g)_{n,m\in  \mz }$ is a frame for $\lr$ with 
frame bounds
$\frac{A}{b},\frac{B}{b}$.
\end{Thm}

Note that if $g$ is compactly supported
and bounded, we get (2) above automatically.  Casazza, Christensen and
Janssen \cite{CCJ} have shown that the conditions in Proposition
\ref{pfcc}
are both necessary and sufficient if $g$ is a real-valued positive
function
on $\mathbb R$.  

A major question in WH-frame theory is:

\begin{Prob}
Classify all the $g,a,b$ so that $(g,a,b)$ is a WH-frame.
\end{Prob}

A special case of this problem which is already extremely difficult
(as we will see in this manuscript) is:

\begin{Prob}
Classify all measurable sets $E\subset \mathbb R$ and $a,b\in \mathbb R$
so
that $({\chi}_{E},a,b)$ is a WH-frame.
\end{Prob}

An even further special case which is still very difficult is  still
open.

\begin{Prob}{\bf [a,b,c-Problem]}
Classify all $a,b,c\in \mathbb R$ so that $({\chi}_{[0,c]},a,b)$ is a
WH-frame.
\end{Prob} 

In these notes we will examine what is known about the latter two
problems and
add some new results to the list.


\section{Compactly Supported Functions}\label{CSF}
\setcounter{equation}{0}

In this section we will look at the more general problem of compactly
supported
functions and WH-frames.  We will not review what is known in this case,
but
instead develop some specific results for use later in examining
characteristic functions which give WH-frames.

Now we give the corresponding result of Proposition \ref{P1}
(5) for two non-zero elements.

\begin{Prop}\label{P2}
Let $g\in \lr$ and assume that for all $0\le t\le 1$ at most two
elements
of $(g(t-n))_{n\in \mz  }$ are non-zero.  The following are equivalent:

(1) $(g,1,1)$ is a WH-frame.

(2)  The CC-condition holds.

(3)  We have $G_{0}\le B< \infty$ a.e. 
and for $H_{n}(t)=|g(t)|-|g(t-n)|$ there is a $0<A$
such that
\[
A \le
\inf \{|H_{n}(t)|: 0\not= n\in \mz  
, g(t)\not= 0 \},
\]
where $(3)$ holds for almost every $t$ with $g(t)\not= 0$. 
\end{Prop}

\begin{proof}
$(1)\Rightarrow (3)$:  By Proposition \ref{P1} (3), we get the
first part of (3).  We will assume the $H_{n}$ condition fails
and show that $(g,a,b)$ fails to have a non-zero lower frame
bound.  In this case,
for a fixed ${\epsilon}>0$, after
a translation, it is
easily seen that there is an $ m\in \mz  $ and a set $E\subset [0,1]$
with $|E|>0$ so that $g(t)\not= 0$ for all $t\in E$ and 
$g(t+m)\not= 0$ for all $t\in E$ and
$$
||g(t)|-|g(t+m)||< \epsilon, \text{ for  all } t\in E.
$$
We construct a function $f\in \lr$ by:
\[
f=f_{0}{\chi}_{E} -f_{1}{\chi}_{E+m}+f_{2}{\chi}_{E+2m}
-f_{3}{\chi}_{E+3m}
+\cdots -f_{2n-1}{\chi}_{E+(2n-1)m}, \]
where $(f_{0},f_{1})$ are given by
$$ 
f_{0}g = |g| \ \ \text{ on } E, \ \ \ 
f_{1}g = |g|\ \ \text{ on } E+m
$$
and $(f_{i})_{i=2}^{(2n-2)m}$ are chosen iteratively
so that there is a 1-periodic function $h_{i}$ with 
$|h_{i}(t)|=1$ and
$$
f\cdot T_{-im}g = h_{i}T_{-im}[|{\chi}_{E}g|-|{\chi}_{E+m}g|].
$$
It follows that for all $i\not= 0, 2n-1$ (since the $h_{i}(t)$ are
1-periodic)
$$
\sum_{k}|<f,E_{k}T_{-im}g>|^{2} =
\sum_{k}|\int_{E}h_{i}(t)[|g(t+m)|-|g(t)|]E_{k}(t)\ dt|^{2}
$$
$$
= \|{\chi}_{E}h_{i}[T_{m}g-g]\|^{2}_{L^{2}(E)} \le {\epsilon}|E|.
$$
 For the other two values of $i$ we will get
$$
\sum_{k}|<f,E_{k}T_{i}g>|^{2}\le \int_{E}|g(t)|\ dt \le |E|B.
$$
So
$$
\sum_{k,n\in \mz  }|<f,E_{k}T_{n}g>|^{2}\le ((n-2){\epsilon}+B)|E|.
$$

Since $\|f\|^{2} = n |E|$, 
it is easily seen that it is now impossible for
$(g,1,1)$ to have a non-zero lower frame bound.

$(3)\Rightarrow (2)$:  We need to check that if $(g,1,1)$ is a 
frame then the conditions of
Theorem \ref{pfcc} are satisfied.  By Proposition \ref{P1}, we know that
for some $A >0$ we have (when $g(t)\not= 0$),
$$
||g(t)|-|g(t+n)||\ge A,
$$
for all $t\in \mathbb R$ and all $0\not= n\in \mathbb Z$.
Hence, 
$$
||g(t)|-|g(t+n)||^{2}\ge A^{2}.
$$
That is,
$$
|g(t)|^{2}+|g(t+n)|^{2}- 2|g(t)||g(t+n)|\ge A .
$$
But, it is easily checked that 
this is precisely the first condition of Theorem \ref{pfcc} for our
case.  Now, if $(g,1,1)$ is a frame then $g$ is bounded and by our
assumptions,
it is easily seen that the second condition must also be satisfied.
\end{proof}

We mention another necessary condition for compactly 
supported functions to give WH-frames.

\begin{Prop}\label{P3}
If $g\in \lr$ has compact support and for every $\epsilon > 0$ there is
a
set $E\subset [0,1]$ with $0< |E|$ and for a.e. $t\in E$ we have
$$
|\sum_{n\in \mz  }g(t+n)|< \epsilon , 
$$
then $(g,1,1)$ is not a WH-frame.
\end{Prop}

\begin{proof}
Without loss of generality we may assume that supp $g\subset [0,N]$.
Fix $n\ge N$ and choose $E\subset [0,1]$ with $|E|>0$ and
$$
|\sum_{k\in \mz  }g(t+k)|\le \frac{1}{n}.
$$
Let
$$
f = \sum_{k=0}^{n}{\chi}_{k+E}.
$$
Then
$$
\|f\|^{2} = n|E|
$$
Also, by our hypotheses, (and mimicking the proof of Proposition
\ref{P2}) for $n-N$ values of $m$ we have
$$
\sum_{k}|<f,E_{k}T_{m}g>|^{2} = |\int_{E}\sum_{k\in \mz  }g(t+k)E_{k}\
dt |^{2}
\le |E|^{2}\frac{1}{n^2}.
$$
For the other terms we get at most:
$$
2(1+2+\cdots (N-1))|E|^{2} = N(N-1)|E|^{2}.
$$
Hence,
$$
\sum_{km\in \mz  }|<f,E_{k}T_{m}g>|^{2}\le
|E|^{2}(N(N-1)+\frac{n-N}{n^{2}})
\le |E|^{2}N^{2}.
$$
So to have a frame, we need an $A>0$ so that
$$
An|E|\le |E|^{2}N^{2}.
$$
That is,
$$
An\le |E|N^{2}\le N^{2}, 
$$
which is a contradiction for large n. 
\end{proof}


\section{the $a,b,c$-Problem}
\setcounter{equation}{0}

  The operator $Lf(t) = f(t/b)$
is an invertible operator on $\lr$ and satisfies:  
$$L(T_{na}{\chi}_{[0,c]})(t)) =
T_{nab}{\chi}_{[o,bc]}(t)
$$ 
while $L({\chi}_{[0,c]}(t) = {\chi}_{[o,bc]}(t)$ and
$L(E_{mb}g) = E_{m}L(g)$.  It follows
that we may as well assume that $b=1$ in the $a,b,c$ Problem.
Some of the results in this section were previously announced in
Casazza \cite{C}.

We will use an immediate consequence of Proposition \ref{P1}
(2).

\begin{Rmk}\label{R}
If $(g,a,1)$ is a WH-frame then $(T_{m}g)_{m\in \mathbb Z}$
is a Riesz basic sequence.
\end{Rmk}

We start with the case $a=b=1$.

\begin{Prop} 
Consider $({\chi}_{[0,c]},1,1)$.  If $c=1$ this is an
orthonormal basis for $\lr$.  If $c<1$ it is a 
normalized tight frame
sequence in $\lr$ which is not a frame.  If $1<c$ it is not a frame
sequence
in $\lr$.
\end{Prop}

\begin{proof}
The case $c=1$ is obvious.  If $c<1$, this is the image
of the orthonormal basis $({\chi}_{[0,1]},1,1)$ under the
orthogonal projection $Pf = {\chi}_{E}f$, where $E =
\cup_{n\in \mathbb Z}([0,c]+n)$.  For $1<c$
we look at two cases:

{\bf Case I}:  $2k-1<c\le 2k$.  Let $d = c-(2k-1)>0$.  For any
natural
number $n>k$ let
$$
f = \sum_{i-0}^{2n}(-1)^{i}{\chi}_{[i,i+d)}.
$$
Now, if ${\ell}\le -(2k)$ or ${\ell}\ge 2n+2k$ then $f,{\chi}_{[0,c]}$
have
disjoint supports so their inner product is 0.  If ${\ell} = 0,1,\cdots
,
2n-2k$ then $<f,{\chi}_{[0,c]}> = 0$.  Otherwise,  $|<f,{\chi}_{[0,c]}>|
=$
 0 or d.  Hence,
$$
\sum_{j}|<f,T_{j}{\chi}_{[0,c]}>|^{2}\le 3(2k)d,
$$
while
$$
\|f\|^{2} = 2nd.
$$
Since $k$ is fixed and $n$ is arbitrary, we have that
$(T_{j}{\chi}_{[0,c]})$ is not a Riesz basic sequence.

{\bf Case II}:  $2k<c\le 2k+1$.  This time let
$$
f =
{\chi}_{[0,d]}-\frac{1}{2}{\chi}_{[1,1+d]}-\frac{1}{2}{\chi}_{[2,2+d]}
+{\chi}_{[3,3+d]}-\frac{1}{2}{\chi}_{[4,4+d]}-\frac{1}{2}{\chi}_{[5,5+d]}
+{\chi}_{[6,6+d]}-\cdots ,
$$
where there are $3n$ terms in the sum above.  Now proceed similarly to 
Case I.    
\end{proof}

We can go further.  

\begin{Rmk}
For the $a,b,c$-Problem, we may assume that $a<b=1<c$.
\end{Rmk}

\begin{proof}
As we saw, we may assume that $b=1$ (and the case
$a=1$ is done).  Hence, by Proposition \ref{P1}
(1) we have that $a\le 1$.  Also, by Proposition
\ref{P1} (4), if $a\le c<1$ then we have a frame and
if $c<a<1$ we do not have a frame.
\end{proof}

Janssen \cite{J1} has shown (using the Walnut representation
of the frame operator) just how delicate the $a,b,c$ Problem is.

\begin{Prop}
Assume $a<1<c$.

(1)  If $a$ is not rational and $1<c<2$ then $({\chi}_{[0,c]},
a,1)$ is a frame.

(2)  If $a=p/q$ is rational, $gcd(p,q)=1$, and $2-(1/q) < c< 2$,
then $({\chi}_{[0,c]},a,1)$ is not a frame.

(3)  If $a>3/4$, $c=L-1+L(1-a)$ with integer $L\ge 3$, then
$({\chi}_{[0,c]},a,1)$ is not a frame.

(4)  If $d$ is the greatest integer $\le c$ and 
$|c-d-1/2|<(1/2)-a$, then $({\chi}_{[0,c]},a,1)$ is a frame.
\end{Prop}

Janssen \cite{J} also has an interesting chart of certain values
where we have a frame or don't have a frame
for several cases of the $a,b,c$-Problem.

We can add to this list the following.

\begin{Prop}
If $2\le c\in \mathbb N$, then for all $a>0$,
$({\chi}_{[0,c]},a,1)$
is not a frame.
\end{Prop}

\begin{proof}
Let $c=n\in \mathbb N$ and $g={\chi}_{[0,n]}$.  Then we have
$$
\|\sum_{j=0}^{k-1}(T_{jn}g-T_{jn+1}g) \|^{2} = 2.
$$
Hence $(T_{n}g)$ is not a Riesz basic sequence and so $(g,a,1)$ cannot
be a WH-frame by Remark \ref{R}. 
\end{proof}


\section{Characteristic Functions Giving WH-Frames}\label{CF1}
\setcounter{equation}{0}

Casazza and Kalton \cite{CK} have shown that the problem of
classifying all characteristic functions which give WH-frames
for the case $a=b=1$ is already exceptionally difficult since it
is equivalent to a classical unsolved problem in complex function
theory.

\begin{Prob}
Classify all integer sets $n_{1}<n_{2}<\cdots n_{k}$ so that
$$
\sum_{j=1}^{k} z^{n_{j}}
$$
does not have any zeroes on the unit circle.
\end{Prob}

For notation, we call  a measurable subset $F$ of $\mathbb R$
{\bf complete} if $|\mathbb R - \cup_{n\in \mathbb Z}(F+n)|=0$, and we
say
that two measurable sets $E,F$ in $\mathbb R$ are {\bf completely
disjoint}
if $|(E+n)\cap (F+m)| = 0$, for all $m,n \in \mathbb Z$.  
We call $E$ a {\bf WH-frame set for a,b} if $({\chi}_{E},a,b)$ is a
frame for $\lr$.  The first
main result of Casazza and Kalton \cite{CK} is:

\begin{Thm}\label{ckth}
Fix integers $n_{1}<n_{2}<\cdots < n_{k}$.  The following are
equivalent:

(1)  The set $F = \cup_{j=1}^{k}([0,1)+n_{j})$ is a Weyl-Heisenberg 
frame set with frame bounds $A,B$.

(2)  We have $A\le  |\sum_{j=1}^{k}z^{n_{j}}|\le B$, for all $|z|=1$,

(3)  For every measurable set $E\subset [0,1]$ of positive measure, and 
$F_{0} = \cup_{j=1}^{k}(E+n_{j})$,  
$(E_{m}T_{n}{\chi}_{F_{0}})_{m,n\in \mathbb Z}$ 
is a frame for $L^{2}(F_{0})$
with frame bounds $A,B$.
\end{Thm}

Casazza and Kalton \cite{CK} also gave an equivalent formulation for
general
characteristic functions to give WH-frames for $a=b=1$.
We call a measurable
set $F\subset \mathbb R$ an {\bf elementary A-Weyl-Heisenberg sub-frame
set of length k} if $F=\cup_{j=1}^{k}(E+n_{j})$ for some $(n_{j})$ and
some measurable subset $E$ in $[0,1)$ and  
$$
A\le \text{inf}_{|z|=1}|\sum_{j=1}^{k}z^{n_{j}}| .
$$

Casazza and Kalton \cite{CK} also classified all WH-frame sets for
$a=b=1$.

\begin{Thm}
Let $F$ be a complete subset of $\mathbb R$.  The following are
equivalent:

(1)  $F$ is a Weyl-Heisenberg frame set.

(2)  There are constants $k,A>0$ so that $F$ can be written as a 
union (finite or infinite) 
of pairwise completely disjoint elementary A-Weyl-Heisenberg sub-frame
sets of length $\le k$.
\end{Thm}

We have the following examples to illustrate the technicalities 
which arise just for the case $a=b=1$.

\begin{Ex}\label{ex02}
If $g = {\chi}_{[0,2]}$ or $g={\chi}_{[0,1]} - {\chi}_{[1,2]}$, then
$(g,1,1)$ is not a frame.  
\end{Ex}
\begin{proof}
This is immediate by Proposition \ref{P2}. 
\end{proof}

\begin{Ex}
If $g = {\chi}_{[0,3]}$, then $(g,1,1)$ is not a frame.
\end{Ex}
\begin{proof}
If we consider the function
$$
f = {\chi}_{[0,1]}-\frac{1}{2}{\chi}_{[1,2]}-\frac{1}{2}{\chi}_{[2,3]}+
{\chi}_{[3,4]}-\frac{1}{2}{\chi}_{[4,5]}-\frac{1}{2}{\chi}_{[5,6]}+{\chi}_{[6,7]}\cdots
,
$$
where the sum has $3n$ terms.  Then $<f,g>=0$ for all terms except 4 of
them
at the ends and these terms yield:  $\frac{1}{2}$ twice and $1$ twice.
So
$$
\sum_{k\in \mz  }|<f,g>|^{2} = \frac{5}{2}.
$$
Since $\|f\| = (n+\frac{n}{2})^{1/2}$, we see that $(g,1,1)$ is not a
frame. 
\end{proof}

\begin{Ex}
If $g = \frac{1}{2}{\chi}_{[0,2]}-{\chi}_{[2,3]}$, then $(g,1,1)$ is
not a frame.
\end{Ex}

\begin{proof}
This is immediate by Proposition \ref{P3}.
\end{proof}

\begin{Ex}
If $g = {\chi}_{[0,2]}-{\chi}_{[2,3]}$, then $(g,1,1)$ does give a 
WH-frame.
\end{Ex}
\begin{proof}
We just note that $(g,1,1)$ satisfies the hypotheses of Theorem
\ref{pfcc}.  In
particular, 
$$
\sum_{n\in \mz  }|g(t-n)|^{2} = 3,
$$
while
$$
G_{k} = 0, \ \ \text{for}\ \ |k|\ge 3,\ \ \text{and}\ \ k=1,-1.
$$
Finally, 
$$
G_{2} = G_{-2} = -1.
$$
So $A = 3-2=1$ and $B = 3+2=5$ in Theorem \ref{pfcc}. 
\end{proof}


\section{Fundamental Frames} 
\setcounter{equation}{0}

 In this section we give some explicit characteristic functions
that yield WH-frames.  Since they come from the very natural
characteristic functions $\chi_{[0,a)}$ and $\chi_{[0,\fb)}$ we refer
to these as the {\bf fundamental frames for the system determined by $a$
and $b$}.  One can use these fundamental frames  to
decompose the frame operator of any PF WH-system $(g,a,b)$.
The main tool of this section  is the
$\fb$-inner product and its norm
\[\label{framerep}\lb f, g\rb_{\fb}(t)=\sum_{k \in \mz }f(t-\fkb
)\overline{g(t-\fkb )}. 
\ \ \ \ \  \|f\|_{\fb}(t)=\sum_{k \in \mz } |f(t-\fkb)|^2 \]

Here we state a few of the necessary properties of this $\fb$-inner
product. First we introduce the standard $\fb$-orthonormal basis. Let
$e_k=T_\fkb\chi_{[0,\fb)}$. .  

\begin{Prop}\label{ip} For all $f,g,h \in \lr$ and $j \in \mz$:

\[ \begin{tabular}{ccccc}
 &(1)&  $\lb f,g \rb_\fb(t) \text{ is } \fb-\text{periodic}$
&(2)&  $ f=\sum_{k\in \mz} \lb f, e_k \rb_\fb(t)  (t) e_k$  \\
&(3)&  $\lb e_0,f\rb_\fb(t) e_0=\ol{f}e_0$  &(4)&  $\lb fg,h\rb_\fb(t)=
\lb
f\overline{h}, \ol{g} \rb_\fb(t)$\\
&(5)&  $\lb T_{\fjb}f,g\rb_\fb(t)= \lb f, T_{-\fjb}g \rb_\fb(t)$
&(6)&  $ T_{\fjb} \lb f,g\rb_\fb(t)= \lb f, g \rb_\fb(t)$ 
  \end{tabular} \] 

\end{Prop}

Ron and Shen \cite{RS1,RS2} and Casazza and Lammers \cite{CL} have
both used inner products of this type to give what the later called
compressions of operators for WH-systems.
Associated with this $\fb$-inner product are a special class of
operators.
\begin{Def}
 We say that a linear
operator $L:\lr \to \lr$ is a $\fb${\bf -factorable
operator} if for any factorization $f={\phi}g$ where $f,g\in
L^{2}(\Bbb R)$ and $\phi$ is an $\fb$-periodic function on $\mr$ 
 we have $$ L(f) =L({\phi}g) = {\phi}L(g).  $$
\end{Def}

We summarize some of the  known results about $\fb$-factorable
operators and compressions in the theorem below. 

\begin{Thm} Let $\|g\|_\fb(t) \le B$ a.e. and  $\|h\|_\fb(t) \le C $a.e.
\\
(1) Let $ L(f) = \sum_{m \in \mz} \lb f, E_{mb} g\rb E_{mb}(t)
h$.  Then $L$ is a $\fb$ factorable operator and has the
following {\bf compression}
\[ \mathcal L(f) =\lb f, g\rb_\fb(t)h\]
(2) The frame operator, frame transform  and preframe  operator 
 for the  system $(g,a,b)$ are $\fb-factorable$
and  may be compressed as follows.
\[S_g(f)=\fb \sum _{n\in \mz }\lb f,T_{na}g \rb_{\fb}(t)T_{na}g, \, \, 
T(f)= \sfb \sum_k\lb f,e_k\rb_\fb(t) T_{ka}g \text { and }
 T^*(f) = \sfb \sum_k\lb f, T_{ka}g\rb_\fb(t) e_k.\]
(3) For any bounded $\fb$-factorable operator $L$ from $\lr$ to $\lr$
\[\|L(f)\|_\fb(t) \le \|L\| \|g\|_\fb(t) \text{ a.e. }.\]
\end{Thm}

In our first application of a frame generated by a characteristic
function  we use the ``frame'' $(\sqrt{b} \chi_{[0,\fb)},\fb,b)$
to produce a pointwise necessary condition for $(g,a,b)$ to be PF. A
moments reflection shows this yields the standard orthonormal basis
for $\lr$ associated with $e^{2 \pi i m b t}$.  However, our concern
will be the connection with the $e_k$'s from above.  \\

The functions $G_k$ have an elementary representation in the
$a$- inner product.  Namely $G_k(t)=\lb g,T_{\fkb}g \rb_a(t)$.
It is well known that if  the system $(g,a,b)$ is PF  then $\sum|G_k|^2
\le B< \infty$. We now interchange
the roles of $a$ and $b$.   

\begin{Thm} If (g,a,b) is PF  with upper frame bound $B$ then 
\[\sum|\lb g,T_{na}g\rb _{\fb}(t)|^2 \le b B\|g\|_{\fb}(t) { a.e.}\]
\end{Thm}  
\begin{proof}  If $(g,a,b)$ is PF we know  by 
Theorem \ref{P1} that $\|g\|_{\fb}(t) \le B$
a.e. and $T$ is a bounded operator.  Hence, because $\|T\|=\|T^*\|$  we
have $\|T^*(g)\|_{\fb}(t) \le \|T\|\|g\|_{\fb}(t)$ and 
\[\lb T^*(g),T^*(g)\rb _{\fb}(t) =
\fb \lb \sum_n \lb g,T_{na}g\rb_\fb e_n, \sum_k \lb g,T_{ka}g\rb_\fb
e_k \rb_\fb=
\fb \sum|\lb g,T_{na}g\rb_{\fb}|^2(t).\]  The last inequality follows
from the $\fb$-orthnormality of the $e_k$.
\end{proof}

If we consider the case $a=b=1$ then this is the same condition as 
$\sum|G_k|^2 \le B<\infty$.  This condition is not sufficient
for having a WH-frame.

\begin{Ex} There exist a system  $(g,1,1)$ such that  $\sum|G_k|^2
\le B<\infty$ yet $(g,1,1)$ is not PF. \end{Ex}

\begin{proof}
 Consider the function
 $g=\sum_{n> 1} \frac{e_n}{n}$. By Casazza, Christensen, and Janssen
\cite{CCJ1} (Corollary 3.7), if $a=b=1$ then a positive real valued
function  
is PF iff   $\sum|G_k| \le B<\infty$ .  However 
a direct computation shows that
$G_k=(1/k)\sum_{n=1}^{k}\frac{1}{n}$ for the above $g$. These $G_k$
are square summable but not summable.  Hence $(g,1,1)$ is not PF.  
\end{proof}

One can present a large number of necessary conditions for the system
(g,a,b) to be PF with these techniques by switching between the frame
operator, preframe  operator and the frame transform. 
We do not know if any of them are also sufficient.   We
present one more representation 
which is stronger than the one above, at least
in the case $a=b=1$.

\begin{Prop} Assume the system (g,a,b) is PF.  Then
\[\lim_{m \to \infty}
\sum_j\frac{|\sum_{k=-m}^{m}T_{ka}g(t-j/b)|^2}{2mb} \le B^2.\] 
\end{Prop}

\begin{proof} Let
$f_m=\sum_{k=-m}^m e_k$ where $e_k$ is the standard $\fb$-orthonormal
basis. Since $T$ is a $\fb$-factorable linear operator from
$\lr$ to $\lr$ we may apply Theorem 5.2 (3) to get:
\[\|T(f_m)\|^2_{\fb} \le B^2 2m \ \ \  { and }\ \ \  
\lb T(f_m),T(f_m)\rb _{\fb}=\fb
\sum_j|\sum_{k=-m}^{m}T_{ka}g(t-j/b)|^2.\] 
Since m was arbitrary we have

\[\lim_{m \to \infty}
\sum_j\frac{|\sum_{k=-m}^{m}T_{ka}g(t-j/b)|^2}{2mb}  \le B^2.\]
\end{proof}

By ``stronger'' we mean that the example presented above ($g=\sum_{n>
1}\frac{e_n}{n}$) does not satisfy
\[\lim_{m \to \infty}
\sum_j\frac{|\sum_{k=-m}^{m}T_{ka}g(t-j/b)|^2}{2mb}  \le B<\infty.\]


\subsection{Frame decompositions}
Now we present a representation of the frame operator for the system
(g,a,b)  that mixes the
a-inner product and the $\fb$-inner product.  
Note that in this theorem we are mixing our a-inner product with the
$\fb$-orthonormal 
basis.  The end result is that each $\lb \lb
T_{-\fjb}g,T_{-\fkb}g\rb _a ,e_0\rb_\fb(t)$
is  a $\fb$ section of an $a$-periodic function which is then extended
$\fb$ periodically.  

\begin{Thm}\label{wrep} Let $(g,a,b)$ be a PF  WH-system.  Then the
frame
operator has the representation 
\bea S(f)=\sum_k\lb f,e_k\rb _{\fb}(t)h_k=\sum_k\lb f,h_k\rb
_{\fb}(t)e_k
 =\sum_k \lb f,e_k\rb _{\fb}(t) \sum_j 
\lb \lb T_{-\fjb}g,T_{-\fkb}g\rb _a,e_0 \rb_\fb(t) e_j\no \eea
where $e_k$ is the standard $\fb$-orthonormal basis and $h_k=S(e_k)$.
  
\end{Thm}

\begin{proof}
Since the system is PF we know that $S$ is a continuous operator from
$\lr $ to $\lr$. Recall $e_k=T_{\fkb} \chi_{[0,\fb)}$ and  
$f=\sum_k\lb f,e_k\rb _{\fb} e_k$.

Now because $S$ is linear, $\fb$-factorable, continuous  and self
adjoint  we get 
\[S(f)=\sum\lb f,e_k\rb _{\fb} S(e_k)=\sum\lb f,S(e_k)\rb _{\fb}e_k\]
So it is enough to look at $S(e_k)$. First let us note that the
computations below rely heavily on the results from
Proposition~\ref{ip}.  Now we compute:
\bea S(e_k) = \sum_j \lb S(e_k),e_j\rb _{\fb} e_j &=&
\sum_j \lb \sum_n \lb e_k,T_{na}g \rb_\fb T_{na}g,e_j\rb_\fb e_j \no
\\
&=&  \sum_j \lb \sum_n \lb e_k,T_{na}g \rb_\fb e_j, \ol{
T_{na}g}\rb_\fb e_j \no \\
&=&  \sum_j \lb \sum_n \lb e_o,T_{na}T_{-\fkb}g \rb_\fb e_0, 
T_{na}T_{-\fjb}\ol{g}\rb_\fb e_j \no \\
&=&  \sum_j \lb \sum_n T_{na}T_{-\fkb}\ol{g}  e_0, 
T_{na}T_{-\fjb}\ol{g}\rb_\fb e_j \no \\
&=&  \sum_j \lb \sum_n T_{na}T_{-\fkb}\ol{g}  
T_{na}T_{-\fjb}{g},e_0\rb_\fb e_j \no \\
&=&  \sum_j \lb \lb T_{-\fjb}g,  
T_{-\fkb}{g}\rb_a,e_0\rb_\fb e_j \no 
\eea

This gives the result.   
 \end{proof}

Now we decompose this frame operator with respect to a natural frame
arising from the parameters $a,b$.  Let
$\alpha_k= T_{ka}\chi_{[0,a]}$. Then by Theorem \ref{P1} (4)and the
representation of the frame operator given in equation (\ref{framerep})
it is clear that $(\sqrt{b}\alpha_0,a,b)$ is a normalized tight frame
and  for all $f\in \lr$ we have 
\[f=\sum_{k\in\mz} \lb f,\alpha_k \rb_\fb(t)\alpha_k.\]  Given the
relation ship between $(g,a,b)$ and $(g,1/b,1/a)$ in 
Theorem \ref{P1} (2) it is also natural to consider the sequence
of characteristic functions $\phi_k=T_\fkb \chi_{[0,a)}$
Putting  this together with the  decomposition above yields the
following.

\begin{Prop}  Let $(g,a,b)$ be a PF WH-system.  The frame operator has
the representation 
\[S(f)=\sum_{k\in \mz} \lb f, \alpha_k\rb_\fb(t)T_{ka}\sum_j G_{-j}(t-\fjb)
\phi_j.\]
\end{Prop}  

\begin{proof} We use the fact that $S$ is $\fb$-factorable and then
apply $S$ to $f=\sum_k \lb f, \alpha_k \rb_\fb(t) \alpha_k$. So we
get  
\[S(f)=\sum_{k \in \mz} \lb f,\alpha_k\rb_\fb(t) S(\alpha_k)\]
It is also well known that $S$ commutes with $T_{ka}$ for all $k\in
\mz$ so it is enough to find $S(\alpha_0)$.  Let us use the
representation above. Since all the $\fb$-inner products are 0 except
one we get:
\[S(\alpha_0)=\lb \alpha_0,e_0\rb_\fb(t)  \sum_j \lb \lb
g,T_{-\fjb}g\rb _a,e_0\rb_\fb e_j .\]  Now since $\lb
\alpha_0,e_0\rb_\fb(t)e_j=\phi_j$  we get the result.  Note that since
the
$\phi_j$ are only supported on an interval of length $a$ we no
longer need the $\fb$-periodic extension.
\end{proof}

\subsection{Equivalent frames and  $S^{\frac{-1}{2}}$}
In this last section we will give one more fundamental example of
WH-frame
obtained
from a characteristic function and another normalized tight frame.
We go on to show that the two frames are equivalent 

First the sake of simplicity we do the case where $\frac{1}{2} \le ab
\le 1$. Let 
$\beta_k=T_{ka}\chi_{[0,\fb)}$.  Then again by  Theorem \ref{P1} (4)
we know that $(\sqrt{b} \beta_0,a,b)$ forms a frame with upper frame
bound 2 and lower frame bound 1 and frame operator
\[S^b(f)=\sum_{k \in \mz}\lb f,\beta_k \rb_\fb \beta_k.\]
  Now let $\delta=\fb-a$ and
$\gamma_0=
\frac{1}{\sqrt{2}}\left(\chi_{[0,\delta)}+
\chi_{[a,\fb)}\right)+\chi_{[\delta,a )}$  
and $\gamma_k=T_{ka} \gamma_0$.  Then  $(\sqrt{b} \gamma_0,a,b)$ forms a
normalized tight frame which yields the following decomposition of the
identity operator:
\[f=\sum_{k \in \mz}\lb f,\gamma_k \rb_\fb \gamma_k.\]

Finally we let $\psi_0=
\frac{1}{\sqrt{2\sqrt{2}}}\left(\chi_{[0,\delta)}+
\chi_{[a,\fb)}\right)+\chi_{[\delta,a )}$, $\psi_k=T_{ka}\psi_0$
  and consider the frame
$(\sqrt{b} \psi,a,b)$ with frame operator
\[S^\psi(f)=\sum_k \lb f, \psi_k\rb_\fb \psi_k.\]

Now we compute  $S^\psi(\beta_0)$
\[S^\psi(\beta_0)=\lb \beta_0, \psi_{-1}\rb_\fb \psi_{-1}+ 
\lb \beta_0, \psi_{0}\rb_\fb \psi_{0}+
 \lb \beta_0, \psi_{1}\rb_\fb \psi_{1}=\gamma_0 \]
Again since frame operators commute with $T_{ka}$ we get that
$S^\psi(\beta_k)= \gamma_k$ and hence the two frames are equivalent.
Furthermore since $\sqrt{b}\gamma_k$ is a normalized tight frame we
have that   
\[S^\psi S^bS^\psi=
S^\psi \left (\sum_{k \in \mz}\lb S^\psi(f),\beta_k
 \rb_\fb \beta_k \right )
=S^\psi \left (\sum_{k \in \mz}\lb f,\gamma_k \rb_\fb \beta_k\right)
=\sum_{k \in \mz}\lb f,\gamma_k \rb_\fb \gamma_k=f
\]
and therefore $S^\psi S^bS^\psi=I$ and hence 
$(S^b)^{\frac{-1}{2}}=S^\psi$

  \end{document}